\documentclass{article}
\usepackage{amsfonts}
\usepackage{epsfig}

\newcommand{\TT}{\mathbb{T}} 
\newcommand{\Li}{\mathrm{Li}} 
\newcommand{\re}{\mathop{\mathrm{Re}}} 
\newcommand{\im}{\mathop{\mathrm{Im}}} 
\newcommand{\ii}{\mathrm{i}} 
\newcommand{\Lf}{\mathrm{L}} 
\newcommand{\HH}{\mathbb{H}} 
\newcommand{\dC}{\mathbb{C}}

\newcommand{\dR}{\mathbb{R}}
\newcommand{\dQ}{\mathbb{Q}}

\newcommand{\dd}{\mathrm{d}}
\newcommand{\e}{\mathrm{e}}
\newcommand{\pf}{\noindent {\bf PROOF.} \quad}
\newcommand{\qed}{$\Box$}
\newtheorem{thm}{Theorem}
\newtheorem{defi}[thm]{Definition}
\newtheorem{prop}[thm]{Proposition}

\newtheorem{ex}[thm]{Example}

\newtheorem{obs}[thm]{Observation}

\begin{document}

\thispagestyle{plain}

\begin{center}
\Large{Mahler Measure and Volumes in Hyperbolic Space}

\medskip
\normalsize{Matilde N. Lal\'{\i}n
{\footnote[1]{E-mail address: mlalin@math.utexas.edu}}
}{\footnote[2]{Supported by Harrington Fellowship}}

\medskip
University of Texas at Austin. Department of Mathematics.
 1 University Station C1200.  Austin, TX 78712, USA
\end{center}

\begin{abstract}The Mahler measure of the polynomials $t(x^m-1) y - (x^n-1) \in \dC[x,y]$ is essentially the sum of volumes of a certain collection of ideal hyperbolic polyhedra in $\HH^3$, which can be determined a priori as a function on the parameter $t$. We obtain a formula that generalizes some previous formulas given by Cassaigne and  Maillot \cite{M} and Vandervelde \cite{V}. These examples seem to be related to the ones studied by Boyd \cite{B1}, \cite{B2} and Boyd and Rodriguez Villegas \cite{BRV2} for some cases of the $A$-polynomial of one-cusped manifolds.

\end{abstract}



\section{Introduction}
The logarithmic Mahler measure of a Laurent polynomial $R \in \dC[x_1^{\pm},\dots, x_n^{\pm}]$ is defined as
\[ m(R) := \frac{1}{(2 \pi \ii)^n} \int_{\TT^n} \log | R(x_1, \dots, x_n)| \frac{\dd x_1}{x_1} \dots \frac{\dd x_n}{x_n} \]
where $\TT^n= \{ (z_1,\dots, z_n)\in \dC^n | |z_1|=\dots = |z_n|=1 \}$ is the $n$-torus.
Jensen's formula provides a simple expression for the Mahler measure of a one-variable polynomial as a function on the roots of the polynomial. The several-variable case is much harder and there are only a few examples of polynomials whose Mahler measure has been found.

The simplest formula for the Mahler measure of a polynomial in more than one variable was found by Smyth \cite{S}:
\begin{equation}
m(x+y+1) = \Lf'(\chi_{-3},-1),
\end{equation}
which expresses the Mahler measure as a special value of the derivative of the L-series in the character $\chi_{-3}$, the nontrivial Dirichlet character of conductor 3.

Later Boyd and Rodriguez Villegas \cite{BRV1} studied the polynomials $R(x,y) = p(x) y - q(x)$. They found that when $p(x)$ and $q(x)$ are cyclotomic, $m(R)$ can be expressed as a sum of values of the Bloch--Wigner dilogarithm at certain algebraic arguments.

The Bloch--Wigner dilogarithm is defined by
\begin{equation}
D(z):= \im(\Li_2(z))+ \log|z|\arg(1-z)
\end{equation}
where $\Li_2(z) = \sum_{n=1}^\infty \frac{z^n}{n^2} $ for $|z| < 1$ is the classical dilogarithm. $D(z)$ can be extended as a real analytic function in $\dC \setminus \{0,1\}$ and continuous in $\dC$.
Let us point out that
\begin{equation}\label{eq:propD}
D(\bar{z}) = - D(z)  \quad ( \Rightarrow D|_{\dR} \equiv 0)
\end{equation}
\begin{equation} \label{eq:prop2D}
-2 \int_0^\theta \log|2 \sin t| \dd t =  D(\e^{2 \ii \theta}) = \sum_{n=1}^\infty \frac{\sin(2n\theta)}{n^2}
\end{equation}
An account of the properties of the Bloch--Wigner dilogarithm can be found in Zagier's work \cite{Z}.

One of these amazing properties is that $D(z)$ is equal to the volume of the ideal hyperbolic tetrahedron of shape $z$, with $\im z >0$ (denoted by $\Delta(z)$). In other words, a tetrahedron in $\HH^3$ whose vertices are $0,1,\infty, z$ (and in particular they belong to  $\partial \HH^3$). See Milnor \cite{Mil}, and Zagier \cite{Za}.

\begin{figure} \label{Maillot}
\centerline{\includegraphics[width=17pc]{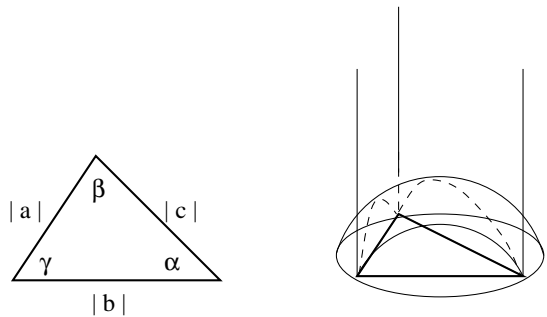}}
\caption{The main term in Cassaigne -- Maillot formula is the volume of the ideal hyperbolic tetrahedron over the triangle.}
\end{figure}

The simplest example of a relation between Mahler measure and dilogarithm (and hence hyperbolic volumes) is given by Cassaigne and Maillot, \cite{M} : for
$a,b,c \in \dC^*$,
\begin{equation}
\pi m(a+bx+cy) = \left \{ \begin{array}{lr} D \left( \left|\frac{a}{b}\right|
\e^{\ii \gamma} \right) + \alpha \log |a| + \beta \log |b| + \gamma \log |c| &
\triangle\\ \\
\pi \log \max \{ |a|, |b|, |c|\} & \mathrm{not}\: \triangle
\end{array}\right.
\end{equation}
where $\triangle$ stands for the statement that $|a|$, $|b|$, and $|c|$ are the
lengths of the sides of a triangle, and $\alpha$, $\beta$, and $\gamma$ are the
angles opposite to the sides of lengths $|a|$, $|b|$, and $|c|$ respectively.
See figure \ref{Maillot}. We see that the dilogarithm term in this formula corresponds to the volume of the ideal tetrahedron that can be built over the triangle of sides $|a|$, $|b|$ and $|c|$.

Another example was considered by Vandervelde \cite{V}. He studied the Mahler measure of $axy+bx+cy+d$ and found a formula, which in the case of $a,b,c,d \in \dR^*$, is very similar to the formula above. The Mahler measure (in the nontrivial case) 
 turns out to be the sum of some logarithmic terms and two values of the dilogarithm, which can be interpreted as the hyperbolic volume of an ideal polyhedron that is built over a cyclic quadrilateral. The quadrilateral has sides of length $|a|$, $|b|$, $|c|$ and $|d|$.

Summarizing,

\begin{itemize}
\item The zero set of Cassaigne -- Maillot 's polynomial is described by
\[ y = \frac{ax + b}{c} \]
and its Mahler measure is the sum of some logarithms and the volume of an ideal polyhedron built over a triangle of sides $|a|$, $|b|$ and $|c|$.
\item The zero set of Vandervelde 's polynomial is described by the rational function
\[ y = \frac{bx + d}{ax + c} \]
and the Mahler measure of the corresponding polynomial is the sum of some logarithms and the volume of an ideal polyhedron built over a quadrilateral of sides $|a|$, $|b|$, $|c|$ and $|d|$.
\end{itemize}

It is natural then to ask what happens in more general cases, for instance, some of the examples given by Boyd and Rodriguez Villegas. We have studied the Mahler measure of
\begin{equation}\label{poly}
R_t(x,y ) = t(x^m - 1) y - (x^n - 1)
\end{equation}
whose zero set is described by the rational function
\begin{equation}
y = \frac{x^{n-1} + \dots + x + 1}{t(x^{m-1}+ \dots + x + 1) }
\end{equation}
We find that the Mahler measure of the polynomial $R_t(x,y ) $ has to do with volumes of ideal polyhedra built over polygons with $n$ sides of length 1 and $m$ sides of length $|t|$. In fact, 
\begin{thm}
\begin{equation}
\pi m(R_t(x,y)) = \pi \log |t| + \frac{2}{mn}\sum \epsilon_k \mathrm{Vol} (\pi^*(P_k)) + \epsilon \sum_{k=1}^N (-1)^k \log |t| \, \arg {\alpha_k}
\end{equation}
where $\epsilon, \epsilon_k = \pm 1$ and the $P_k$ are  all the admissible polygons of type $(m,n)$.
\end{thm}

Here $R_t$ is the polynomial (\ref{poly}) and admissible polygons are, roughly speaking, all the possible cyclic polygons that can be built with $n$ sides of length 1 and $m$ sides of length $|t|$. We will give a precise definition later.

This formula with hyperbolic volumes is similar to certain formulas that occur for some cases of the $A$-polynomial of one-cusped manifolds. This situation was studied by Boyd \cite{B1}, \cite{B2}, and Boyd and Rodriguez Villegas \cite{BRV2}.

We have divided this paper as follows. In section 2 we compute the Mahler measure of the family of polynomials $t(x^m-1)y -(x^n-1)$. In section 3, we relate this formula to volumes in hyperbolic space. In section 4, we show some examples of this relationship and in section 5 we explore analogies to the $A$-polynomial situation.

\section{A preliminary formula}
\begin{prop} \label{prop} Consider the polynomial
\[ R_t(x,y) = t(x^m-1) y - (x^n-1), \qquad t \in \dC^* \qquad gcd(m,n)=1\]
Let $\alpha_1, \dots \alpha_N \in \dC$ be the different roots (with odd
multiplicity) of
\[ Q(x) = \frac{x^n-1}{x^m-1} \cdot \frac{x^{-n}-1}{x^{-m}-1} - |t|^2 \]
such that $|\alpha_k|=1$, $\alpha_k \in \HH^2 = \{ z \in \dC\, | \im z >0\}$  , and they are ordered
counterclockwise starting from the one that is closest to 1. Then
\begin{equation} \label{eq:prop1}
\pi m(R_t(x,y)) = \pi \log |t| + \epsilon \sum_{k=1}^N (-1)^k \left( \frac {D({\alpha_k}^n)}{n}-
\frac {D({\alpha_k}^m)}{m}  + \log |t| \, \arg {\alpha_k} \right)
\end{equation}
where $\epsilon = \pm 1$.
\end{prop}

\pf This Proposition is very similar to Proposition 1 in \cite{BRV1} (when $t=1$), but we prove it here so we can provide more details. We may suppose
that $t \in \dR_{>0}$, since multiplication of $y$ by
numbers of absolute value 1 does not affect the Mahler measure (an easy property that can be deduced from the definition of Mahler measure).  By Jensen's formula,

\begin{eqnarray}
2 \pi m(R_t(x,y))- 2\pi \log t & = & \frac{1}{\ii} \int_{\TT^1} \log^+ \left| \frac
{1-x^n}{t(1-x^m)}
\right| \frac {\dd  x}{x}\\
& = & \frac{1}{\ii} \sum_j \int_{\gamma_j}  \log \left|
\frac {1-x^n}{t(1-x^m)} \right| \frac {\dd  x}{x} \label{eq:prop}
\end{eqnarray}

\begin{figure}
\centerline{\includegraphics[width=15pc]{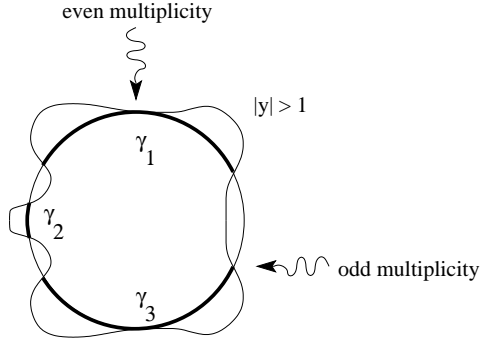}}
\caption{The arcs $\gamma_i$ are the sets where $|y| \geq 1$. The
  extremes of these arcs occur in points where $y$ crosses the
  unit circle. }
\label{f:prop}
\end{figure}

Where $\log^+ x = \log x$ for $x>1$ and 0 otherwise. Here $\gamma_j$ are the arcs of the unit circle where $\left| \frac
{1-x^n}{t(1-x^m)} \right| \geq 1$. The extreme points of the
$\gamma_j$ must be roots of
$Q(x)$. It is easy to see that we only need to consider the roots of odd
multiplicity, indeed, $y= \frac{1-x^n}{t(1-x^m)} $ crosses the
unit circle only on those roots. See figure \ref{f:prop}.

It is also clear that for each  root of $Q(x)$, its inverse is
also a root (in other words, $Q$ is reciprocal), so the roots with absolute value one come in conjugated pairs, except, maybe, for $1$ and
$-1$. We need to analyze what happens with these two cases.

Case $-1$. Since $m$ and $n$ are coprime, they cannot be even at the same time
and the only meaningful case in which $-1$ may be a root is when both are odd.
In that case,  $Q(-1) = 1 - t^2$, and $t=1$. Studying the multiplicity
of $-1$ in this case is equivalent to studying the multiplicity of $-1$ as a root of $Q_1(x)
= x^m+x^{-m} - x^n - x^{-n}$. It is easy to see that $Q_1'(-1) =0$ and
$Q_1''(-1) \not = 0$, hence $-1$ is a root of multiplicity two.

Case $1$. We have $Q(1) = \frac{n^2}{m^2} -t^2$. Hence $1$ is root of $Q(x)$
if and only if $t = \frac{n}{m}$.
As before, it is enough to study the parity of the multiplicity of $1$ as a
root of $Q_1(x) = m^2(x^n -2 + x^{-n}) - n^2(x^m -2 + x^{-m})$.  Again we
see that $Q_1'(1) = Q_1''(1) = Q_1'''(1) =0$ but $Q_1^{(4)}(1) \not = 0$.
Hence $1$ is a root of even multiplicity.

Thus we do not need to take $1$ or $-1$ into account and the extremes of the
$\gamma_j$ will  lie in the $\alpha_k$ and their conjugates.

We have
\[\int_{\alpha_k}^{\alpha_{k+1}} \log | 1 - x^n| \frac {\dd  x}{ \ii
  x}  =\frac{1}{n}\int_{\alpha_k^n}^{\alpha_{k+1}^n} \log | 1 - y| \frac {\dd  y}{ \ii
  y}  =  \frac{D({\alpha_k}^n) - D({\alpha_{k+1}}^n)}{n}\]
by equation(\ref{eq:prop2D}).
Using this in formula (\ref{eq:prop}), and with the previous observations about the
roots of $Q$, we obtain formula (\ref{eq:prop1}).
\qed

There is another way of performing this computation, which was suggested by Rodriguez Villegas. The idea is to start with the case of $t=1$ and obtain the general case as a deformation.

In order to do this, let us compute the initial case is a slightly different way. Recall that
\[R_1(x,y) = (x^m-1)y-(x^n-1) \]
Following \cite{RV}, let $X$ be the smooth projective completion of the complex zero locus of $R_1$ and let $S$ the set of points where either $x$ or $y$ have a zero or a pole. Consider the following differential 1-form in $X \setminus S$:
\begin{equation}
\eta(x,y) := \log |x| \dd \arg y - \log |y| \dd \arg x
\end{equation}
Observe that
\begin{equation}
2\pi m(R_1(x,y)) = \int_{\gamma_1} \eta(x,y)
\end{equation}
where $\gamma_1 = \cup \gamma_{1,i}$ is, as before, the set in the unit circle where $|y| \geq 1$. Because of Jensen's formula, the identity above is true for any polynomial in two variables such that the "coefficient" of the highest power of $y$ (a polynomial in $x$) has Mahler measure zero (i.e., product of cyclotomic polynomials and powers of $x$).

It is easy to check that $\eta$ is closed. When the polynomial is tempered (see \cite{RV}) $\eta$ extends to the whole $X$. If in addition $\{ x, y\} = 0$ in $K_2(X)\otimes \dQ$, then the form is exact and in fact,
\begin{equation}
\eta(x, 1-x) = \dd D(x)
\end{equation}
The above conditions are easily verified by $R_1$. Using the above equation, we recover the statement for $t=1$:
\[\int_{\alpha_k}^{\alpha_{k+1}} \log | 1 - x^n| \frac {\dd  x}{ \ii
  x} =   -\frac{1}{n} \int_{\alpha_k}^{\alpha_{k+1}}\eta(x^n,1-x^n) =  \frac{D({\alpha_k}^n) - D({\alpha_{k+1}}^n)}{n}\]

\begin{figure}
\centerline{\includegraphics[width=10pc]{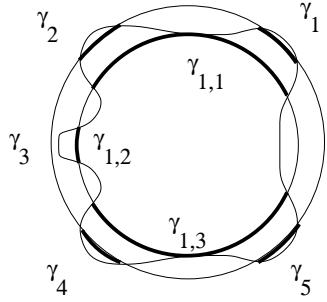}}
\caption{While the arcs $\gamma_i$ are the sets where $|y| \geq 1$, the arcs $\gamma_i$ are the sets where $|y| \geq t$.}
\label{f:fer}
\end{figure}

Now in order to treat the general case, we write
\[ y' = \frac{y}{t} \]
and
\[R_t(x,y') = R_1(x,y)\]
Then if
\[ 2\pi m(R_1) = \int_{\gamma_1} \eta(x,y) \]
where $\gamma_1 = \cup \gamma_{1,i}$,  we also have
\[ 2 \pi m(R_t) = 2\pi \log t + \int_{\gamma} \eta(x,y') \]
where $\gamma = \cup \gamma_j$. Now $\gamma$ is the set of the circle where $|y'| \geq 1$, i.e., where $|y| \geq t$. See figure \ref{f:fer}.

Observe that
\[ \eta(x,y') = \eta(x,y) - \eta(x,t) = \eta(x,y) + \log t \, \dd \arg x \]
Then

\[2 \pi m(R_t) = 2 \pi \log t + \int_{\gamma} \eta(x,y) + \int_{\gamma} \log t \, \dd \arg x \]
\[= 2 \pi \log t + 2 \, \epsilon  \sum_{k=1}^N (-1)^k \left( \frac {D({\alpha_k}^n)}{n}-
\frac {D({\alpha_k}^m)}{m}  + \log t \, \arg {\alpha_k} \right)\] \qed

The moral of this last procedure is that we can compute the Mahler measure of the  polynomial with general coefficient $t$ by altering the Mahler measure of the polynomial with $t=1$.
This concept of smooth deformation will appear in the main result of this paper, when we interpret this formula as volumes in hyperbolic space.

We see that in order to compute the general Mahler measure, we need to integrate in a different path, i.e., in the set where $|y| \geq t$. There is no reason to think that this integration is harder to perform than the one over the set $|y| \geq 1$. However, determining the new arcs $\gamma$ as functions of $t$ might be hard.  Nevertheless, this method should be useful to compute other examples, but we will not go on into this direction in this paper.

In general, it seems difficult to interpret the intersection points (i.e., the starting and ending points for the arcs $\gamma$) geometrically. The main point of this paper is to show such an interpretation for the particular example  that we have studied.

\section{The main result}

We will need some notation. The following definition is not standard.

\begin{defi} A cyclic plane polygon $P$ will be called admissible of type $(m,n)$ if the following conditions are true:
\begin{itemize}
\item $P$ has $m+n$ sides, $m$ of length $t$ (with $t \in \dR_{>0}$) and $n$ of length 1.
\item All the sides of length $t$ wind around the center of the circle in the same direction, (say counterclockwise), and all the sides of length 1 wind around the center of the circle in the same direction, which may be opposite from the direction of the sides of length $t$ (so they all wind counterclockwise or clockwise).
\end{itemize}
\end{defi}

\begin{figure}
\centerline{\includegraphics[width=17pc]{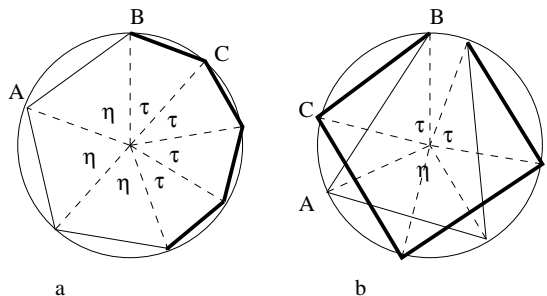}}
\caption{From now on, the bold segments indicate sides of length $t$,
  which are opposite to angles measuring $\tau$. The ordinary segments
  indicate sides of length $1$, opposite to angles measuring $\eta$. The circles in the pictures may seem to have the same radius, but that is not true. The radius is determined by the size polygon that is inscribed in the circle.}
\label{thm1}
\end{figure}

In order to build such a polygon $P$, we need to define two angles, $\eta$ and
$\tau$, which are the central angles subtended by the chords of
lengths $1$ and $t$ respectively. See figure \ref{thm1}. The polygon
does not need to be convex or to wind exactly once around the center
of the circle. Figure \ref{thm1}.a shows an ordinary convex polygon
winding once. Observe that in this picture the two families of sides wind in the same direction.
In figure \ref{thm1}.b, the polygon does not wind around the center. In this picture the two families of sides wind in opposite directions.

Let us remark that there are finitely many admissible polygons for given $m$, $n$ and $t$. Given a polygon, the radius of the circle is fixed. Conversely, for each radius, there is at most one admissible polygon of type $(m,n)$ that can be inscribed in the circle for $t$ fixed. It is easy to see that the radius $r$ and the parameter $t$ satisfy an algebraic equation. So for given $t$ there are only finitely many solutions $r$.

\begin{figure}
\centerline{\includegraphics[width=17pc]{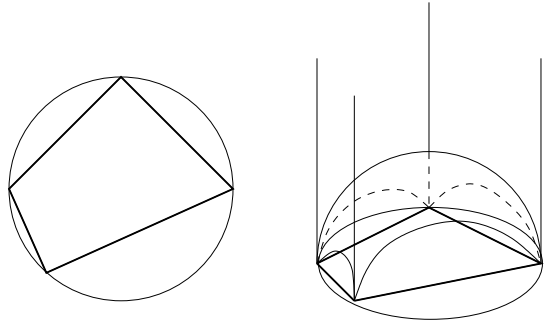}}
\caption{This picture shows an example of how to build the ideal polyhedron
over an admissible polygon.}
\label{f:thm}
\end{figure}

Given an admissible polygon $P$, we think of $P \subset \dC \times \{ 0\} \subset  S^2_\infty  \cup
\HH^3$, then $\pi^*(P)$ denotes the ideal polyhedron whose vertices are  $\infty$ and those of $P$
(see figure \ref{f:thm}). We use the model $\HH^3  \cong \dC \times \dR_{>0} \cup \{\infty \}$ for the hyperbolic space.

Now, it makes sense to speak of the hyperbolic volume $\mathrm{Vol}(\pi^*(P))$ (up to a sign).  We consider $P$ to be subdivided into $m+n$ triangles, as can be seen in figure \ref{thm1}, all of the triangles sharing a vertex at the center of the circle, and the opposite side to this vertex being one of the $m+n$ chords. Hence we get $m$ isosceles triangles of basis $t$ and $n$ isosceles triangles of basis $1$. We consider the orthoschemes over each of these triangles and the total volume will be the sum of the volumes of these tetrahedra, but we take the tetrahedra over the sides of length $t$ to be negatively oriented if the two families of sides wind in opposite directions. Compare with the definition of dilogarithm of an oriented cyclic quadrilateral given by Vandervelde, \cite{V}.

We are now ready to state our main result.

\begin{thm} \label{theo} The dilogarithm term of the Mahler measure in
formula (\ref{eq:prop1}) is equal to the sum of the volumes of certain ideal polyhedra
in the hyperbolic space $\HH^3$:
\begin{equation} \label{eq:thm}
\epsilon \sum_{k=1}^N (-1)^k \left( \frac {D({\alpha_k}^n)}{n}-
\frac {D({\alpha_k}^m)}{m}\right)= \frac{2}{mn}\sum \epsilon_k \mathrm{Vol} (\pi^*(P_k))
\end{equation}
where $\epsilon_k = \pm 1$ and the $P_k$ are  all the admissible polygons of type $(m,n)$.
\end{thm}

\pf First, we will see that for each $\alpha = \alpha_k$, there exists an admissible polygon $P$ as in the statement such that
\begin{equation} \label{eq:vol}
\pm\, 2 \mathrm{Vol} ( \pi^*(P)) =  m D({\alpha}^n) - n D({\alpha}^m)
\end{equation}

Suppose $\alpha = \e^{\ii \sigma}$ is such that
\begin{equation}
\sigma \in  \left(\frac{k \pi}{m} , \frac {(k+1) \pi}{m} \right]   \bigcap  \left(\frac{l \pi}{n} ,
\frac{(l+1) \pi}{n} \right]  \qquad 0 \leq k < m, \quad  0 \leq l < n
\end{equation}

We choose $\eta$ and $\tau$ according to the rules given by the following tables:

\bigskip

\renewcommand{\arraystretch}{1.27}
\begin{tabular}{cc}
\begin{tabular}{|l|l|}
\hline
$k$ even & $\eta := m \sigma - k \pi$ \\
$k$ odd  &  $\eta := (k+1) \pi - m \sigma$ \\
\hline
\end{tabular}&
\begin{tabular}{|l|l|}
\hline
$l$ even & $\tau := n \sigma - l \pi$ \\
$l$ odd  &  $\tau := (l+1) \pi - n \sigma$ \\
\hline
\end{tabular}
\end{tabular}
\renewcommand{\arraystretch}{}

\bigskip
\noindent This choice of $\eta$, $\tau$ is the only possible one that  satisfies
\[\begin{array}{ccc}
\eta & \equiv & \pm\, m \,\sigma \,\mathrm{mod}\, 2\pi \\
\tau & \equiv & \pm\, n \, \sigma \, \mathrm{mod}\, 2 \pi
\end{array}\]
in addition to $0 < \eta \leq \pi$ and $0 < \tau \leq \pi$. The above congruences will guarantee the right arguments in the dilogarithm. This will be clearer later.

Also note that we cannot have $\eta=\tau=\pi$, since this would imply that
\[\frac{k+1}{m} = \frac{l+1}{n}\]
and this is possible only when $\sigma=\pi$, since $m$ and $n$ are coprime and $k<m$, $l<n$. But we have already seen that $\sigma < \pi$ in the Proof of Proposition \ref{prop}.

 Let us prove that such a polygon with these angles and sides does exist. We have that $n \eta \pm m \tau = h \, 2 \pi$. Then, in order that the polygon can be inscribed in a circle, it is enough to verify the Sine Theorem. Take the triangle $\stackrel{\triangle}{ABC}$ in figure
\ref{thm1}. The side $\overline{AB}$, of length 1, is opposite to an angle measuring
$\frac{\eta}{2}$ or $ \pi - \frac{\eta}{2}$. The side $\overline{BC}$ has length
$t$ and is opposite to an angle measuring $\frac{\tau}{2}$ or $ \pi -
\frac{\tau}{2}$. By the Sine Theorem,

\begin{equation}\label{eq:sin}
\frac{1}{\sin \frac{\eta}{2}} = \frac{t}{ \sin \frac{\tau}{2}}
\end{equation}
Looking at the table above, equality (\ref{eq:sin}) becomes
\begin{equation}
\frac{1}{\left| \sin\frac{m \sigma}{2}\right |} = \frac{t}{ \left| \sin \frac{n\sigma}{2} \right |}
\end{equation}
Squaring and using that $1-\cos \omega= 2 \sin^2 \frac{\omega}{2}$,
\begin{equation}
2-2\cos n \sigma = t^2(2- 2\cos m \sigma)
\end{equation}
Since $\alpha = \e^{\ii \sigma}$, we get
\begin{equation}
\left| \frac{\alpha^n -1}{\alpha^m -1}\right| = t
\end{equation}
Since this equation is the algebraic relation $Q(\alpha) = 0$ satisfied by
$\alpha$.

Hence, given $\sigma$, we can find $\eta$ and $\tau$, and
we are able to construct the polygon $P$. The equality $n \eta \pm m \tau = h\, 2
\pi$ indicates that the polygon winds $h$ times around the center of the
circle. The possible sign "--" should be interpreted as a change in the
direction we are going in the circle (from clockwise to counterclockwise or
vice versa, as explained in the definition of admissible polygons). Note that the admissible
polygon constructed in this way is unique, up to the rigid transformations on the plane and
up to the order we choose for the angles $\eta$ and $\tau$ as we wind around
the center. There is no need to place all the $\eta$ first and then all
the $\tau$. We do that for the sake of simplicity and coherence in the
pictures.

Let us prove that of the volume of the
corresponding hyperbolic object is given by formula (\ref{eq:vol}). As we mentioned above, the polygon is divided by $m+n$ triangles, all of them sharing one vertex at the
center of the circle. The volume of the orthoscheme over each of these triangles depends
only on the central angle $\omega$ and is equal to $\frac{D(\e^{\ii \omega})}{2}$,
according to Lemma 2, from Appendix, in Milnor's work \cite{Mil}. We get the following

\renewcommand{\arraystretch}{1.5}
\begin{center}
\begin{tabular}{cc}
\begin{tabular}{|l|r|}
\hline
$k$ \mbox{even}, $l$ \mbox{even} & $\pm ( n D(\e^{\ii \eta}) - m D(\e^{\ii \tau}))$  \\\hline
$k$ \mbox{even}, $l$ \mbox{odd} & $\pm ( n D(\e^{\ii \eta}) + m D(\e^{\ii \tau}))$ \\
\hline
$k$ \mbox{odd}, $l$ \mbox{even} & $\pm ( n D(\e^{\ii \eta}) + m D(\e^{\ii \tau}))$ \\
\hline
$k$ \mbox{odd}, $l$ \mbox{odd} &  $\pm( n D(\e^{\ii \eta}) - m D(\e^{\ii \tau}))$ \\
\hline
\end{tabular} & $= \pm (n D(\alpha^m) - m D(\alpha^n))$\\
\end{tabular}
\end{center}
\renewcommand{\arraystretch}{}

Now we will study the converse problem: given an admissible polygon $P$, determined by
$\eta$ and $\tau$ (i.e., with $n$ sides opposite to the angle $\eta$ and $m$ sides opposite to the angle $\tau$, and $gcd (m,n) =1$). We want to find the $\alpha$ that corresponds to $P$.  The equation describing the polygon $P$ is
either $n \eta + m\tau = h\, 2  \pi$ or $n \eta - m\tau = h\, 2 \pi$. Without loss of generality, we may suppose that $\eta < \tau \leq \pi$ (and so, $0<\eta< \pi$).

Let $s$ be such that
\begin{equation}
s \, n \equiv h \, \mathrm{mod}\, m
\end{equation}
where $s$ is chosen uniquely in such a way that
\begin{equation}
0 <  \eta - s\, 2 \pi < m\, 2 \pi
\end{equation}
The condition $0 < \eta < \pi$ guarantees that $\eta -s \, 2 \pi
\not = m \, \pi$. There are two cases:
\[0 < \eta - s\, 2 \pi < m\, \pi \quad \Rightarrow \quad \sigma := \frac{\eta -
s\, 2 \pi}{m}\]
\[m\, \pi <\eta - s\, 2 \pi < m\, 2 \pi \quad \Rightarrow \quad \sigma := \frac{(s+m)\, 2 \pi - \eta}{m}\]
It is easy to see that these choices work, in the sense that $Q(\alpha)=0$, $\alpha \in \HH^2$ and  $|\alpha| =1$.  We get $0 < \sigma <
\pi$ in both cases, and $\eta \equiv \pm \, m \sigma (\mbox{mod}\, 2
\pi)$, the sign being the one that we need to have the inequality $0 < \eta \leq \pi$, so
we recover the $\alpha$ that produces $\eta$ according to the table above we used to construct $\eta$. What happens with $\tau$? We have $n \eta \pm m \tau = h\, 2 \pi$, then
\[ \sigma = \frac{\eta - s\, 2 \pi}{m}\quad \Rightarrow \quad \pm
\tau =  - n \sigma + \frac{h-sn}{m}\, 2 \pi\]
\[ \sigma = \frac{(s+m)\, 2 \pi - \eta}{m} \quad \Rightarrow \quad
\pm \tau = n \sigma + \left( \frac{h-sn}{m} -n\right ) 2 \pi \]

So that $\tau \equiv \pm\, n \sigma (\mbox{mod}\, 2 \pi)$, and everything is consistent.

\qed

Let us observe that Theorem \ref{theo} gives another proof for the finiteness of the number of admissible polygons (since they all correspond to roots of $Q$).

\section{Some examples}

Let us do some examples illustrating the situation of Theorem \ref{theo}.
\begin{ex} Consider the case $m=2$, $n=3$:
\begin{equation}
y = \frac{x^3-1}{t (x^2-1)}
\end{equation}
\end{ex}

The advantage of this particular case is that we can compute the
actual values of $\alpha_k$. In fact, clearing (cyclotomic) common factors,
\[ R_t(x,y) = t(x+1) y - (x^2+x+1)\]
We need  $|y|=1$, this is equivalent to
\[\frac{x^2+x+1}{t(x+1)} \cdot \frac{x^{-2}+x^{-1}+1}{t(x^{-1}+1)} = 1 \]
The $\alpha_k$ are the roots (in $\HH^2$, with absolute value 1) of the above equation, which can be
expressed as
\[Q_1(x) = x^4 + (2- t^2) x^3 +(3-2t^2)x^2 +(2-t^2) x + 1 = 0 \]

As always, we can suppose that $t>0$.
We know that $Q_1$ is reciprocal. Because of that, we may write, (by simple inspection),
\begin{equation}
Q_1(x) = x^2 S(x+x^{-1}), \quad \mbox{where} \quad  S(M) = M^2+(2-t^2)M + 1-2t^2
\end{equation}

The roots of $Q_1$ in the unit circle in $\HH^2$ correspond to roots of $S$ in the interval $[-2,2]$ and vice versa. This is because  if $\beta = \e^{\ii \theta}$ is a root of $Q_1$, it corresponds to $M= \beta + \beta^{-1} = 2 \cos \theta$, a root of $S$.

The roots of $S$ are
\begin{equation}\label{eq:M}
M = \frac{t^2 - 2 \pm t \sqrt{t^2 +4}}{2}
\end{equation}

We see that $M \in \dR$ always. We also have the following:
\begin{eqnarray}
\left|  \frac{t^2 - 2 - t \sqrt{t^2 +4}}{2} \right| \leq 2 & & \forall
\; 0<t \nonumber\\
\left|  \frac{t^2 - 2 + t \sqrt{t^2 +4}}{2} \right| \leq 2 & & \forall
\; 0<t \leq \frac{3}{2} \nonumber \\
\end{eqnarray}

We obtain either one or two pairs of roots of the form $\{ \alpha, \bar{\alpha} \}$
according to the number of solutions for $M$. Indeed, $\re \alpha =
\frac{M}{2}$ if $M \in [-2,2]$.

The situation is summarized by:
\begin{obs} Let $\alpha_1, \alpha_2$ be such that $\re \alpha_1 =  \frac{t^2 - 2 - t \sqrt{t^2 +4}}{4} $  for $0 < t $,  $\re \alpha_2 = \frac{t^2 - 2 + t \sqrt{t^2 +4}}{4}  $ for  $0<t < \frac{3}{2}$, $|\alpha_i| =1 $ and $\im \alpha_i >0$ . Then for $\sigma_i = \arg \alpha_i$, we have

\begin{equation}\label{eq:sigma1}
\pi > \sigma_1 > \frac{2 \pi}{3}
\end{equation}

\begin{equation} \label{eq:sigma2}
\frac{2 \pi}{3} > \sigma_2 > 0
\end{equation}

\end{obs}

\begin{figure}
\centerline{\includegraphics[width=7.2pc]{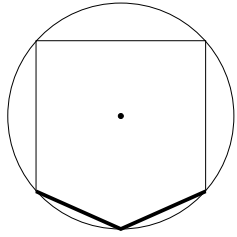}}
\caption{The case of $\alpha_1$ corresponds to the ordinary convex
  polygon. Note that $\alpha_1$ exists  for any $t >0$ and the same is
  true for the polygon.}
\label{ex2}
\end{figure}

We will apply the procedure given in the proof of Theorem \ref{theo} in order to get the
polygons. The case of $\alpha_1$ is very simple.  Because of inequality
(\ref{eq:sigma1}),  $k=1$ and $l=2$ always, so $\eta = 2 \pi - 2
\sigma_1$ and $\tau = 3 \sigma_1 - 2 \pi$. Then $3\eta + 2 \tau = 2 \pi $. This
corresponds to the convex pentagon which is inscribed  in a circle (see figure \ref{ex2}).

The case of $\alpha_2$ splits into three subcases according to the
values of $t$, as shown in the following table.

\begin{table}[h]
\caption{Case $\alpha_2$ }\label{t1}
\renewcommand{\arraystretch}{1.5}
\setlength\tabcolsep{1.5mm}
\begin{tabular}{c|c|l|l|l}
\hline
$0 < t < \frac{1}{\sqrt{2}}$ & $\frac{2 \pi}{3} > \sigma_2 > \frac{\pi}{2} $ &
$k = 1$, $l = 1$ & \begin{tabular}{ccl} $\eta$ & = & $2 \pi -2 \sigma_2$ \\ $\tau$
& = & $2 \pi - 3 \sigma_2$ \end{tabular} & $3 \eta - 2 \tau = 2 \pi$\\
\hline
$\frac{1}{\sqrt{2}} < t < \frac{2}{\sqrt{3}}$ & $ \frac{\pi}{2} > \sigma_2 >
\frac{\pi}{3}$ & $k = 0 $, $l = 1$
& \begin{tabular}{ccl} $\eta$ & = & $2 \sigma_2$
\\ $\tau$ & = & $2 \pi - 3 \sigma_2$ \end{tabular} & $3\eta + 2\tau = 4 \pi$ \\
\hline
$\frac{2}{\sqrt{3}} < t < \frac{3}{2}$ & $ \frac{\pi}{3} > \sigma_2 > 0$ &
$k = 0$, $l = 0$
 & \begin{tabular}{ccl} $\eta$ & = & $2 \sigma_2$ \\ $\tau$ & = &
$3 \sigma_2$ \end{tabular} & $3\eta - 2\tau = 0$ \\
\hline
\end{tabular}
\renewcommand{\arraystretch}{}
\end{table}

\begin{figure}
\centerline{\includegraphics[width=26pc]{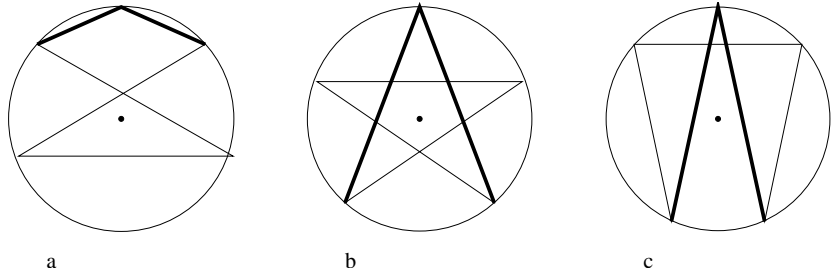}}
\caption{Case $\alpha_2$: a) $\quad 0 < t < \frac{1}{\sqrt{2}}\quad $
  b) $\quad \frac{1}{\sqrt{2}} < t < \frac{2}{\sqrt{3}}\quad $ c)
  $\quad \frac{2}{\sqrt{3}} < t < \frac{3}{2}\quad$ }
\label{ex1}
\end{figure}

Figure \ref{ex1} illustrates the polygons corresponding to each of
these subcases.

We would like to point out that in every case, the condition over $t$ that
we use to compute $\eta$ and $\tau$ is the same condition that assures
that we can build the corresponding polygon.

\begin{figure}
\centerline{\includegraphics[width=17pc]{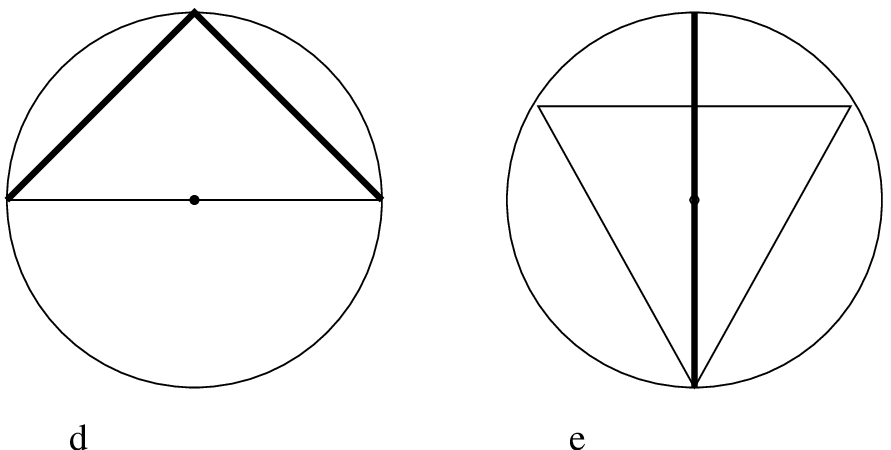}}
\caption{d)$\quad t = \frac{1}{\sqrt{2}}\quad $ e) $\quad t = \frac{2}{\sqrt{3}}\quad$}
\label{exlim}
\end{figure}

The cases of $t = \frac{1}{\sqrt{2}}$ and $t = \frac{2}{\sqrt{3}}$ are
limit cases and we get the transition pictures of figure
\ref{exlim}.

Note that figure \ref{exlim}.d is indeed the intermediate figure
between \ref{ex1}.a and \ref{ex1}.b and the same is true for
\ref{exlim}.e, which is between \ref{ex1}.b and \ref{ex1}.c.

\qed

The next example is harder to examine.

\begin{ex}Consider the case $m=1$, $n=4$:
\begin{equation}
y = \frac{x^4-1}{t (x-1)}
\end{equation}
\end{ex}
which corresponds to
\[R_t(x,y) = t y -(x^3+x^2+x+1) \]
The difficulty in this case lies in the fact that we are not able to compute the values of the $\alpha_k$ exactly.

We proceed as in the other example. We need to study the roots in $\HH^2$, with absolute value one, of the polynomial:

\[ Q_1 (x) =x^6 + 2 x^5 + 3 x^4+(4-t^2)x^3 + 3 x^2 + 2 x + 1 = 0 \]

As always, we can suppose that $t>0$. As before, $Q_1$ is reciprocal. Because of that, we may write
\begin{equation}
Q_1(x) = x^3 S(x+x^{-1}), \quad \mbox{where} \quad  S(M) = M^3 + 2 M^2 - t^2
\end{equation}

As before, the roots of $Q_1$ in the unit circle in $\HH^2$ correspond to roots of $S$ in the interval $[-2,2]$ and vice versa.

Solving this cubic equation would lead to too many calculations. Instead, we will analyze in what cases we get $M \in \dR$ with $|M|\leq 2$. First of all, observe that
\[S'(M) = 3 M^2 + 4 M\]
Then $S$ decreases in $\left (- \frac{4}{3}, 0 \right ) $ and increases everywhere else. Since $S(0) < 0$, it  has always one positive real root. For this root to be in $(0,2)$ we need $S(2) > 0$, i.e., $t< 4$. Hence for $0 < t < 4$ we get $0 < M < 2$ and
\begin{equation}
 \frac{\pi}{2} > \sigma_1 > 0
\end{equation}

On the other hand, $S$ has two negative roots iff $S\left( - \frac{4}{3} \right ) > 0$, i.e.,  $t^2 < \frac{32}{27}$. Since $S(-2) < 0$, both roots are in $(-2,0)$. Hence for $0 < t < \sqrt{ \frac{32}{27}}$, we get two roots, $-2 < M < - \frac{4}{3}$ and $ - \frac{4}{3} < M < 0$ corresponding to
\begin{equation}
\pi > \sigma_2 > \arccos \left( -\frac{2}{3} \right)
\end{equation}
and
\begin{equation}
 \arccos \left( -\frac{2}{3}\right)  > \sigma_3 > \frac{\pi}{2}
\end{equation}
 respectively.

Let us build the polygons corresponding to the three solutions $\alpha_1$, $\alpha_2$, $\alpha_3$. When $t > 4$ the Mahler measure is $\log t$ (the terms with $D$ vanish) and this is consistent with the fact that we do not get any polygons.

Since $m=1$, $k=0$ always and $\eta = \sigma_i$.

Consider the case of $\alpha_1$. It splits into two subcases according to the values of $t$.

\begin{table}[h]
\caption{Case $\alpha_1$}
\renewcommand{\arraystretch}{2}
\begin{center}
\begin{tabular}{c|c|l|l|l}
\hline
$0 < t < \sqrt{4 + 2 \sqrt{2}} $& $\frac{\pi}{2} > \sigma_1 > \frac{\pi}{4} $ &
 $l=1$ & $\tau =  2 \pi - 4 \sigma_1$ & $4 \eta +  \tau = 2 \pi$\\
\hline
$  \sqrt{4 + 2 \sqrt{2}} < t < 4 $ & $ \frac{\pi}{4} > \sigma_1 > 0 $ & $l=0$ & $\tau =  4 \sigma_1$& $4 \eta - \tau = 0 $ \\
\hline
\end{tabular}
\end{center}
\renewcommand{\arraystretch}{}
\end{table}

\begin{figure}
\centerline{\includegraphics[width=26pc]{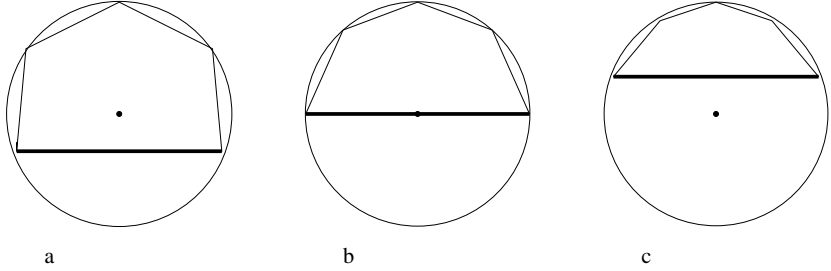}}
\caption{Case $\alpha_1$: a) $\quad 0 < t < \sqrt{4+2\sqrt{2}} \quad $
  b) $\quad t = \sqrt{4+2\sqrt{2}} \quad $ c)
  $\quad \sqrt{4+2\sqrt{2}} < t < 4 \quad$ }
\label{polygon_exb1}
\end{figure}

These two subcases correspond to an ordinary polygon. See figure \ref{polygon_exb1}. When $\sigma_1 = \frac{\pi}{4}$ we get a transition between the two pictures, when the side of length $t$ coincides with the diameter of the circle (figure  \ref{polygon_exb1}.b). In order to compute the value of $t$, we consider the usual subdivision of the polygon into  triangles, all of them  sharing a vertex at the center of the circle. This subdivision involves five triangles, one for each side of the polygon, but we obtain just four triangles in this case, because one triangle degenerates in the diameter of the circle. Each of these four triangles has a side of length 1, opposite to an angle $\sigma_1 = \frac{\pi}{4}$, and the other two sides of length $\frac{t}{2}$. By the Sine Theorem,
\[ \frac{t}{2 \sin \frac{3 \pi}{8}} = \frac{1}{\sin \frac{\pi}{4}} = \sqrt{2} \]
on the other hand,
\[ \frac{1}{\sqrt{2}}= \sin \frac{\pi}{4} = \sin \frac{3 \pi}{4} = 2 \sin \frac{3 \pi}{8} \cos \frac{3 \pi}{8} \]
then,
\[ \frac{1}{\sqrt{2}}= 2 \frac{t}{2 \sqrt{2}} \sqrt{1 - \frac{t^2}{8}} \]
It is easy to solve the above equation. We need a solution that is positive and such that $\frac{t}{2} > 1$ (because larger sides are opposite to larger angles in a triangle). Such a solution corresponds to \[t = \sqrt{4 + 2 \sqrt{2}}\]

Now let $0 < t < \sqrt{ \frac{32}{27}}$. For $\alpha_2$, we have two subcases.

\begin{table}[h]
\caption{Case $\alpha_2$}
\renewcommand{\arraystretch}{2}
\setlength\tabcolsep{1.15mm}
\begin{center}
\begin{tabular}{c|c|l|l|l}
\hline
$0 < t < \sqrt{4 - 2 \sqrt{2}} $& $\pi > \sigma_2 > \frac{3\pi}{4} $ &
 $l=3$ & $\tau =  4 \pi - 4 \sigma_2$ & $4 \eta +  \tau = 4 \pi$\\
\hline
$  \sqrt{4 - 2 \sqrt{2}} < t < \sqrt {\frac{32}{27}} $ & $ \frac{3\pi}{4} >
\sigma_2 > \arccos\left( -\frac{2}{3}\right) $ & $l=2$ & $\tau =  4 \sigma_2 - 2 \pi$& $4 \eta - \tau = 2 \pi $ \\
\hline
\end{tabular}
\end{center}
\renewcommand{\arraystretch}{}
\end{table}

\begin{figure}
\centerline{\includegraphics[width=26pc]{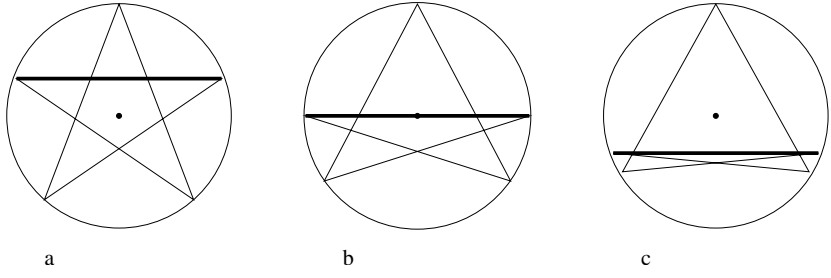}}
\caption{Case $\alpha_2$: a) $\quad 0 < t < \sqrt{4-2\sqrt{2}} \quad $
  b) $\quad t = \sqrt{4-2\sqrt{2}} \quad $ c)
  $\quad \sqrt{4-2\sqrt{2}} < t < \sqrt{\frac{32}{27}} \quad$}
\label{polygon_exb2}
\end{figure}

These subcases correspond to stars. When $\sigma_2 = \frac{3 \pi}{4}$ we get the transition picture, when the side of length $t$ coincides with the diameter of the circle. (Figure \ref{polygon_exb2}.b). In order to compute the value of $t$, we subdivide the picture into four triangles as we did with the transition case for $\alpha_1$. Each of these triangles has a side of length 1, opposite to an angle $\sigma_2= \frac{3\pi}{4}$, and the other two sides of length $\frac{t}{2}$. By the Sine Theorem,
\[ \frac{t}{2 \sin \frac{ \pi}{8}} = \frac{1}{\sin \frac{3\pi}{4}} = \sqrt{2} \]
using that
\[ \frac{1}{\sqrt{2}}= \sin \frac{\pi}{4} = 2 \sin \frac{ \pi}{8} \cos \frac{ \pi}{8} \]
we get the same equation as before
\[ \frac{1}{\sqrt{2}}= 2 \frac{t}{2 \sqrt{2}} \sqrt{1 - \frac{t^2}{8}} \]
Now we take into account that the solution must be positive and must satisfy that $\frac{t}{2} < 1$, so it must be \[t = \sqrt{4 - 2 \sqrt{2}}\]

For $\alpha_3$, we have one case with two kinds of pictures.

\begin{figure}
\centerline{\includegraphics[width=26pc]{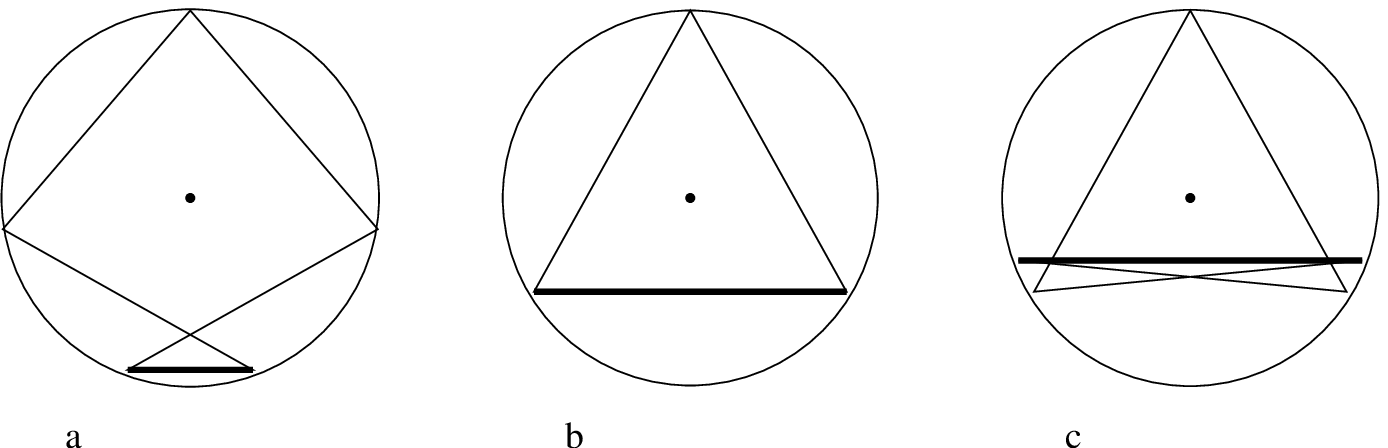}}
\caption{Case $\alpha_3$: a) $\quad 0 < t < 1\quad $
  b)$\quad t = 1 \quad $ c)
  $ \quad 1  < t < \sqrt{ \frac{32}{27}} \quad$ }
\label{polygon_exb3}
\end{figure}

\begin{table}[h]
\caption{Case $\alpha_3$}
\renewcommand{\arraystretch}{2}
\begin{center}
\begin{tabular}{c|c|l|l|l}
\hline
$0 < t < 1 $& $ \frac{\pi}{2} < \sigma_3 < \frac{2\pi}{3}$ & $l=2$ & $\tau =   4 \sigma_3 -2 \pi $ & $4 \eta -  \tau = 2 \pi$\\
\hline
$1 < t < \sqrt{\frac{32}{27}}$& $\frac{2\pi}{3} <  \sigma_3 < \arccos \left( - \frac{2}{3}\right)  $ &
 $l=2$ & $\tau =  4 \sigma_3 - 2 \pi $& $4 \eta - \tau =  2 \pi $ \\
\hline
\end{tabular}
\end{center}
\renewcommand{\arraystretch}{}
\end{table}
We get two possible pictures, see figure \ref{polygon_exb3}. The transition picture is a triangle, when $t=1$ and two of the other sides coincide with the side of length $t$ (see figure \ref{polygon_exb3}.b).

The second subcase corresponds to a star as in the second subcase of $\alpha_2$. If we look carefully, we will notice that the angle $\tau$ is different in both cases. The geometric interpretation for this fact is that there are two ways of building the star for $\sqrt{4 - 2\sqrt{2}} < t < \sqrt{\frac{32}{24}}$.

\qed

\begin{ex} Consider the general case with $t=1$:
\begin{equation}
y = \frac{x^n-1}{x^m-1}
\end{equation}
\end{ex}
This is one particular case of the polynomials studied in \cite{BRV1}.

Without loss of generality, we can suppose $n > m$ (since the Mahler measure remains invariant under the transformation $y \rightarrow y^{-1}$). We need to look at
\[ Q(x) = x^n(x^n - x^m - x^{-m} + x^{-n})  = (x^{m+n} - 1) ( x^{n-m} -1)\]
It is easy to see that the roots of $Q$ are  $\zeta_{m+n}$ and $\zeta_{n-m}$, the $m+m$ and $n-m$ roots of the unity.

Getting the pictures is a delicate task, involving considerations such as the parity of $m$ and $n$.
We will content ourselves with studying the case $m=1$. Then the roots of $Q$ are $\zeta_{n+1}$ and $\zeta_{n-1}$. We only need the roots in $\HH^2$.  In other words:
\[ \sigma_j = \frac{2 j \pi }{n+1} \quad \mbox{for} \quad j=1, \dots , \left[ \frac{n}{2} \right] \]
\[ \rho_j = \frac{2 j \pi }{n-1} \quad \mbox{for} \quad j=1, \dots , \left[ \frac{n}{2} \right] -1 \]

Since $m=1$, $k=0$ always. The choice for $\eta$ and $\tau$ is given by the following table:

\begin{table}[h]
\caption{$\eta$ and $\tau$}
\renewcommand{\arraystretch}{1.7}
\begin{center}
\begin{tabular}{c|l|c|l}
\hline
$ \sigma_j =  \frac{2 j \pi}{n+1} $ &$ l = 2 j - 1$ & \begin{tabular}{ccl} $\eta$ & = & $\sigma_j$ \\ $\tau$
& = & $2 j \pi  - n \sigma_j$ \end{tabular} ($\eta = \tau =  \frac{2 j \pi}{n+1} $)& $n \eta + \tau = 2 j \pi  $ \\
\hline
$\rho_j = \frac{2 j \pi}{n-1} $ & $ l = 2 j $&  \begin{tabular}{ccl} $\eta$ & = & $\rho_j$ \\ $\tau$
& = & $n \rho_j - 2 j \pi $ \end{tabular} ($\eta = \tau =  \frac{2 j \pi}{n-1} $) & $ n \eta - \tau = 2 j \pi  $ \\
\hline
\end{tabular}
\end{center}
\renewcommand{\arraystretch}{}
\end{table}

We can see that we get different kinds of stars, and that all the sides wind around the center of the circle in the same direction for the $\sigma_j$ and that the two families of sides wind in different directions for the $\rho_j$. For this simple relation to hold, it is crucial that $m=1$. The general case is much harder to describe.

\qed

\section{Analogies with the case of $A$-polynomials}

Boyd \cite{B1}, \cite{B2} and Boyd and Rodriguez Villegas \cite{BRV2} found several examples where the Mahler measure of the $A$-polynomial of a compact, orientable, complete, one-cusped, hyperbolic manifold $M$ is related to the volume of the manifold. The $A$-polynomial is an invariant $A(x,y) \in \dQ[x,x^{-1}, y, y^{-1}]$.
Boyd and Rodriguez Villegas found identities of the kind
\[ \pi m(A) = \mathrm{Vol} (M) \]

Motivated by those works, we wonder if there is any relation with our situation. Consider the case of $t=1$. Then the terms with $\log |t|$ vanish and formulas (\ref{eq:prop1}) and (\ref{eq:thm}) become
\begin{equation} \label{eq:A}
\pi m(R_1(x,y)) = \frac{2}{mn} \sum \epsilon_k \mathrm{Vol} ( \pi^*(P_k))
\end{equation}

Let us first mention a few words about $A$-polynomials. The $A$-polynomial is a certain invariant from the space of representations $ \rho : \pi_1(M) \rightarrow SL_2(\dC)$, more precisely, it is the minimal, nontrivial algebraic relation between two parameters $x$ and $y$ which have to do with $\rho(\lambda)$ and $\rho(\mu)$, where $\lambda, \mu \in \pi_1(\partial M)$ are the  longitude and the meridian of the boundary torus. For details about this definition see for instance, \cite{CCGLS}, \cite{CL1}, \cite{CL2}.

Assume the manifold $M$ can be decomposed as a finite union of ideal tetrahedra:
\begin{equation}
M = \bigcup_{j=1}^k \Delta(z_j)
\end{equation}
Then
\begin{equation}
\mathrm{Vol} (M) = \sum_{j=1}^k D(z_j)
\end{equation}

For this collection of tetrahedra to be a triangulation of $M$, their parameters must satisfy certain equations, which may be classified in two groups:
\begin{itemize}
\item {\em Gluing equations\/}. These reflect the fact that the tetrahedra all fit well around each edge  of the triangulation:
\begin{equation}\label{eq:glue}
\prod_{i=1}^k z_i^{r_{j,i}} (1 -z_i)^{r_{j,i}'} = \pm 1 \qquad \mbox{for} \quad j=1,\dots, k
\end{equation}
where $r_{j,i}$, $r_{j,i}'$ are some integers depending on $M$.
\item{\em Completeness equations\/}. These have to do with the triangulation fitting properly at the cusps. If there is one cusp, there will be two of them:

\begin{eqnarray}
\prod_{i=1}^k z_i^{l_i} (1 -z_i)^{l_i'} & = & \pm 1\\
\prod_{i=1}^k z_i^{m_i} (1 -z_i)^{m_i'} & = & \pm 1
\end{eqnarray}

where $l_i$, $l_i'$, $m_i$, $m_i'$ are some integers depending on $M$.
\end{itemize}

One possible solution to this system of equations is the geometric solution, when all the $\im z_i > 0$. There are other possible solutions. For the geometric solution, $\sum_{i=1}^k D(z_i)$ is the volume of $M$. Following Boyd \cite{B2}, $\sum_{i=1}^k D(z_i)$ will be called a pseudovolume of $M$ for the other solutions. Hence, pseudovolumes correspond to sums where at least one of the terms is the dilogarithm of a number $z$ with $\im(z) \leq 0$. One way to understand this is that some of the tetrahedra may be degenerate (when $\im z=0$) or negatively oriented (when $\im z <0$). Boyd shows some examples where
\[ \pi m(A) = \sum V_i \]
where $V_0 = \mathrm{Vol} (M)$ and the other $V_i$ are pseudovolumes. At this point, the analogy of our situation with Boyd's results should be clear. We would like to say that $R_1$ is some sort of $A$-polynomial for some hyperbolic object.

Back to the construction of the $A$-polynomial, introduce "deformation parameters" $x$ and $y$ and  replace the completeness relations by

\begin{eqnarray}
\prod_{i=1}^k z_i^{l_i} (1 -z_i)^{l_i'} & = & x^2 \label{eq:l}\\
\prod_{i=1}^k z_i^{m_i} (1 -z_i)^{m_i'} & = & y^2 \label{eq:m}
\end{eqnarray}

For our purposes, the $A$-polynomial is obtained by eliminating $z_1,\dots, z_k$ from the system formed by the gluing equations (\ref{eq:glue}) and the equations (\ref{eq:l}) and (\ref{eq:m}). This construction is slightly different form the original definition, since it parameterizes representations in $PSL_2(\dC)$ instead of $SL_2(\dC)$. For a detailed discussion about the relationship between this definition and the original one, we refer the reader to Champanerkar's thesis \cite{Ch}. See also Dunfield's Appendix to \cite{BRV2}.

Following \cite{NZ}, we form the matrix
\begin{equation}
U= \left( \begin{array}{cccccc} l_1 & \cdots & l_k & l_1'& \cdots & l_k' \\ m_1 & \cdots & m_k & m_1' & \cdots & m_k'\\ r_{1,1} & \cdots & r_{1,k} & r_{1,1}' & \cdots & r_{1,k}'\\ \vdots &\ddots &\vdots & \vdots & \ddots & \vdots \\r_{k,1} & \cdots & r_{k,k} & r_{k,1}' & \cdots & r_{k,k}'
\end{array} \right)
\end{equation}

One of the main results in \cite{NZ} is
\begin{thm}(Neumann--Zagier)\label{thm:NZ}
\begin{equation}
U J_{2k} U^t = 2 \left( \begin{array}{cc}J_2 & 0 \\ 0 & 0 \end{array} \right)
\end{equation}
where
\[J_{2p} = \left( \begin{array}{cc} 0 & I_p \\ -I_p & 0 \end{array} \right) \]
\end{thm}

Back to our problem, we recall that each term in the right side of equation (\ref{eq:A}) corresponds to the volume of an orthoscheme that is built over an admissible polygon. Each of these polygons is naturally divided into $m+n$ triangles. This division yields a division of the corresponding orthoscheme in $m+n$ hyperbolic tetrahedra. These tetrahedra are not ideal. However, we can redo the whole process by pushing this common vertex that lies over the center of the circle to the base plane $\dC \times \{ 0 \}$ and all the tetrahedra become ideal. The new orthoscheme will be denoted by $\pi^*(P_k')$. The volumes of the tetrahedra get multiplied by 2. Formula (\ref{eq:A}) becomes

\begin{equation} \label{eq:B}
\pi m(R_1(x,y)) = \frac{1}{mn} \sum \epsilon_k \mathrm{Vol} ( \pi^*(P'_k))
\end{equation}

Inspired by the above situation, it is natural for us to take these tetrahedra as a triangulation for our hyperbolic object. So we would like to choose the shape parameters to be $w=\e^{\ii \eta}$ and $z = \e^{\ii \tau}$. Here we actually mean that we have $k=m+n$ tetrahedra, $m$ of them have parameter $w$ and $n$ of them have parameter $z$. We choose the parameters to be $w_1, \dots, w_m$ and $z_1, \dots ,z_n$ and impose the additional condition that $w_1 = \dots = w_m$ and $z_1= \dots =z_n$. The fact that the tetrahedra wind around the axis through the center of the circle which is orthogonal to the base plane $\dC \times \{ 0 \}$, can be expressed by the gluing equation $w_1 \dots w_m z_1 \dots z_n =1$. Further, we need two additional completeness equations, which will be chosen ad hoc for the final result to fit our needs.

It is easy to see that the system
\begin{equation}\label{eq:tri}
\left \{
\begin{array} {r c l} w_1^{\alpha} z_1^{\beta} &  =  & x^2 \\
w_1^{-mn(m+n)\alpha} z_1^{-mn(m+n)\beta}(1-w_1)^{2n} \dots (1-w_m)^{2n} & & \\
 \cdot (1-z_1)^{-2m} \dots (1-z_n)^{-2m}  & = & y^2 \\
w_1 \dots w_m z_1 \dots z_n & = & 1 \\
w_1 w_2^{-1}& = & 1\\
\vdots \\
w_1w_m^{-1}&  = &1\\
z_1 z_2^{-1} & = & 1\\
\vdots \\
z_1z_n^{-1} & = &1\\
\end{array} \right .
\end{equation}
for $n \alpha - m \beta=1$, satisfies the conditions of Theorem \ref{thm:NZ}.  The system reduces easily to

\begin{equation}\label{eq:tri2}
\left \{
\begin{array} {rcl} w^{\alpha} z^{\beta}  & = & x^2 \\
w^{-mn(m+n)\alpha} z^{-mn(m+n)\beta}\left( \frac{1-w}{1-z}\right) ^{2mn} & = & y^2 \\
w^mz^n & = & 1 \\
\end{array} \right .
\end{equation}

Replace the first equation by its $(n^2-m^2)$th-power,

\begin{equation}\label{eq:tri3}
\left \{
\begin{array} {rcl} w^n z^m   & = &  x^{2(n^2-m^2)} \\
\left(\frac{z}{w} \right)^{mn} \left(\frac{1-w}{1-z}\right)^{2mn} & =&  y^2 \\
w^mz^n & = & 1 \\
\end{array} \right .
\end{equation}

Eliminate $x$ and $y$. One of the branches (the one with $w=x^{2n}, z=x^{-2m}$), is
\begin{equation}
y^2 = \left(\frac{x^n-x^{-n}}{x^m-x^{-m}}\right)^{2mn}
\end{equation}
Consider
\begin{equation}
\tilde{R}(x,y) = (x^m-x^{-m})^{mn} y - (x^n-x^{-n})^{mn}
\end{equation}
(We have chosen a particular branch again). It is easy to see that
\[mn\cdot m(R_1) = m(\tilde{R}) \]
Hence
\begin{equation}
\pi m(\tilde{R}(x,y)) = \sum \epsilon_k \mathrm{Vol} ( \pi^*(P'_k))
\end{equation}

We can think of $\tilde{R}$ as the $A$-polynomial of some hyperbolic object that has a triangulation that can be described by the system of equations (\ref{eq:tri}).  We do not expect this object to be a manifold. For instance, this object cannot be the complement of a knot, since the $A$-polynomial of a knot has a number of properties such as being reciprocal (\cite{CL1}, \cite{CL2}).

Also note that the objects whose volumes we are adding, correspond to solutions of the system of equations (\ref{eq:tri3}) with $x=y=1$, in other words, we are able to recover the $\alpha_k$ for the case of $t=1$.
In fact, we need to solve the system
\begin{equation}
\left \{
\begin{array} {rcl} w^n z^m   & = & 1 \\
\left(\frac{z}{w} \right)^{mn} \left(\frac{1-w}{1-z}\right)^{2mn} & =&  1 \\
w^mz^n & = & 1 \\
\end{array} \right.
\end{equation}
From the third equation we write $w = u^n$, $z= u^{-m}$. Substituting this into the first equation, we see that $u^{n^2-m^2} =1$. From this we conclude that $|w|=|z|=1$. Now look at the second equation, which says that
\[ \left( \frac{(1-w)(1-w^{-1})}{(1-z)(1-z^{-1})}\right )^{mn} = 1\]
If we take into account that $|w|=|z|=1$, we see that we are actually computing the $mn$ -- power of an absolute value, then
\[ \frac{(1-w)(1-w^{-1})}{(1-z)(1-z^{-1})} =1 \]
Concluding that the only possible solutions are $m+n$ and $|n-m|$ - roots of unity is now an easy exercise.

Note that we generally get more than one "geometric solution", in the sense that there is more than one solution where all the parameters lie in $\HH^2$. For instance, in the case of $(m,n)=(2,3)$, we get the two solutions described in the example. In the case of $(m,n) =(1,4)$, we get two "geometric solutions" and one that corresponds to $\alpha_3$ with $t=1$ (the triangle in figure \ref{polygon_exb3}.b) which is not "geometric" since $\tau$ must be taken to wind in the opposite direction (which is equivalent to $\im z < 0$).

Let us also observe that, if $X$ is the smooth projective completion of the curve defined by $\tilde{R}(x,y)=0$, then
$\{ x, y \} = 0$ in $K_2(X) \otimes \dQ$. In fact, in $\bigwedge^2 ( \dC(X)^*) \otimes \dQ$,
\[ x \wedge y = mn \, x \wedge (x^n-x^{-n})  - mn \, x \wedge (x^m-x^{-m})  \]
\[= m x^n \wedge (x^{2n}-1) + n x^{-m} \wedge (1-x^{-2m}) \]
Now use that $w = x^{2n}$ and $z=x^{-2m}$ up to torsion,
\[ =\frac{m}{2} w \wedge (1-w) + \frac{n}{2} z \wedge (1-z) \]

The identity above reflects the concept of triangulation as it appears in \cite{BRV2}.

Finally, we should  point out that we could have done the whole process starting from the family
of polynomials $ \tilde{R}_t(x,y) = t (x^m-x^{-m})^{mn} y - (x^n-x^{-n})^{mn}$. The reason we started with the polynomials $R_t(x,y) = t (x^m - 1)y - (x^n - 1)$ is that they are easier to analyze and that they seem a natural choice in order to generalize previous works as it is explained in the introduction.
\bigskip

\section{Final Remark}

To conclude, let us mention that it would be interesting to generalize these examples to cases involving $k$ parameters $t_1, \dots, t_k$ instead of just one parameter. In these general examples we would have to consider admisible polygons with some sides of length $t_1$, some sides of length $t_2$ and so on, thus generalizing completely the examples considered by Cassaigne -- Maillot and Vandervelde.

\bigskip
\noindent{\bf Acknowledgements}
Thanks are due to Fernando Rodriguez Villegas for several helpful discussions and for the idea that started this project. I am grateful to him as well as David Boyd  for their support and encouragement. Finally, I would like to thank the Referee whose comments have led to several improvements in the exposition of this paper.

\end{document}